\documentclass[11pt]{article}
\usepackage{amsmath,amssymb,amsfonts}
\newcommand{\R}{\mathbb{R}}

\newcommand{\fg}{\mathfrak g}
\newcommand{\fh}{\mathfrak h}
\newcommand{\fk}{\mathfrak k}

\newcommand{\fn}{\mathfrak n}

\newtheorem{definition}{Definition}[section]
\newtheorem{lemma}[definition]{Lemma}
\newtheorem{theorem}[definition]{Theorem}

\newtheorem{cor}[definition]{Corollary}

\newtheorem{pro}[definition]{Proposition}

\def\co {cohomogeneity one }

\begin{document}

\vspace*{5mm}

\begin{center} {\bf \Large Cohomogeneity One Dynamics on Three Dimensional Minkowski Space}

{\bf P. Ahmadi  }\\
{ \it Department of Mathematics, University of Zanjan\\
P.O.Box 45371-38791, Zanjan,Iran }\\
email:p.hmadi@znu.ac.ir

\end{center}

\vspace{5mm}

{\small \leftskip 2cm \rightskip 2cm \baselineskip 5mm

{\bf Abstract.}
In this paper we give a classification of closed and connected Lie
groups, up to conjugacy in $Iso({\R^3_1})$, acting by \co on the three
dimensional Minkowski space $\R^3_1$ in both cases, proper and nonproper dynamics. Then we determine causal properties and types of the orbits.

\vspace{5mm} {\bf Keywords}:  {\it Cohomogeneity one , Minkowski space.}

\hspace{18mm} 2010 {\it Mathematics Subject Classification}:
53C30, 57S25.}

\section{Introduction}
There is a basic question involved most types of modern dynamics: which groups admit actions of the type under investigation. A number of researchers have tried to give a list of groups admitting nonproper (or orbit nonproper) actions on Lorentz manifolds, up to local isomorphism or at most up to isomorphism (see for example
\cite{Adams1, Kow, Zeg, Zim}). Many mathematicians have studied \co Riemannian manifolds (see \cite{Br, MK, PT, PS, S}). Most of the works of the second group have concentrated on the study of geometrical properties of the results of the action, such as, existence of slices, the orbits up to isometry, the orbit space up to homeomorphism, etc. In these works the common hypothesis is that the acting group is a closed and connected Lie subgroup of $Iso_g(M)$, where $g$ denotes the Riemannian metric on the smooth manifold $M$. This hypothesis cause an strong dynamical restriction, that is the action should be proper. When the metric is indefinite, this assumption in general
does not imply that the action is proper, so the study becomes much
more difficult. Also some of the results and techniques of the
definite metric fails for the indefinite metric. In this paper, we have attempted to give the list of closed and connected Lie groups, up to conjugacy in $Iso(\R^3_1)=O(1,2)\ltimes \R^3$, acting isometrically and by \co on the three dimensional Minkowski space $\R^3_1$. Then we have tried to determine, for each group in the list, wether its action is proper or nonproper. Finally we have studied causal properties and types of the orbits, for both cases, proper and nonproper actions. Also, we have specified the orbit space, when the action is proper. In the case that the action is nonproper, the orbit space may not be Hausdorff, and so the study seems to be not interesting. Furthermore, we could not use the same definition of principal and singular orbit which is used in \cite{AA}. So we have used a new definition which is compatible with that of proper actions (see preliminary section).

As an interesting result of this paper is that if the action is proper, then the linear part of the acting group is either compact or hyperbolic one parameter subgroup, i.e. it has no nilpotent element. Another considerable result is about the existence of exceptional orbits. It is a well known result that, for proper actions, each exceptional orbit is
nonorientable if the $G$-manifold $M$ is orientable and if $M$ is
simply connected, there is no exceptional orbit (see \cite [p.185]
{Br}), but we see in Propositions \ref{pro-nonproper-1} to \ref{pro-nonproper-3}
that for $M=\R^3_1$, which is simply connected and orientable, there
are orientable exceptional orbits!

\section{preliminaries}

Let $G$ be a Lie group which acts on a connected smooth manifold $M$. The Lie algebra of $G$ is denoted by $\fg$. For each point $x$ in $M$, $G(x)$ denotes the orbit of $x$, and $G_x$ is the stabilizer in $G$ of $x$. The manifold $M$ is called of {\it \co} under the action of the Lie
group $G$ if an orbit has codimension one. The action is said to be {\it proper} if the mapping
$\varphi:G\times M\rightarrow M\times M , \ (g,x)\mapsto(g.x,x)$ is
proper. Equivalently, for any sequences $x_n$ in $M$ and $g_n$ in $G$, $g_nx_n\rightarrow y$ and $x_n\rightarrow x$ imply that $g_n$ has a convergent subsequence. The $G$-action on $M$ is {\it nonproper} if it is not proper. Equivalently, there are sequences $g_n$ in $G$ and $x_n$ in $M$ such that $x_n$ and $g_nx_n$ converge in $M$ and $g_n\rightarrow \infty$, i.e. $g_n$ leaves compact subsets. For instance, if $G$ is compact, the action is obviously proper. The action of $G$ on $M$ is proper if and only if there is a complete
$G$-invariant Riemannian metric on $M$ (see \cite {Alee}). This theorem makes a link between proper actions and Riemannian $G$-manifolds. The orbit
space $M/G$ of a proper action of $G$ on $M$ is Hausdorff, the
orbits are closed submanifolds, and the stabilizers are
compact (see \cite{Adams}). The orbits $G(x)$ and $G(y)$ have the {\it same orbit type} if
$G_x$ and $G_y$ are conjugate in $G$. This defines an equivalence
relation among the orbits of $G$ on $M$. Denote by $[G(x)]$ the
corresponding equivalence class, which is called the {\it orbit type} of
$G(x)$. A submanifold $S$ of $M$ is called a slice at $x$ if there is a $G$-invariant open neighborhood $U$ of $G(x)$ and a smooth equivariant retraction $r:U\rightarrow G(x)$, such that $S=r^{-1}(x)$. A fundamental feature of proper actions is the
existence of slice (see \cite{PT2}), which enables one to define a partial ordering
on the set of orbit types. The partial ordering on the set of orbit types is
defined by, $[G(y)] \leq [G(x)]$ if and only if $G_x$ is conjugate
in $G$ to some subgroup of $G_y$. If $S$ is a slice at $y$, it
implies that $[G(y)] \leq [G(x)]$ for all $x\in S$. Since $M/G$
is connected, there is a largest orbit type in the space of orbit types. Each
representative of this largest orbit type is called a principal
orbit. In other words, an orbit $G(x)$ is principal if and only if
for each point $y\in M$ the stabilizer $G_x$ is conjugate to
some subgroup of $G_y$ in $G$. Other orbits are called singular.
We say that $x\in M$ is a principal point if $G(x)$ is a principal
orbit.

But for the nonproper action there is not slice in general, so we
can not use the same definitions required the existence of slices as
before, hence we use the definition 2.8.1 of \cite{DK} for
determining the principal, singular or exceptional orbits. According
to it for the action of a Lie group $G$ on a smooth manifold $M$,
The points $x,y \in M$, are said to be of the same type, with
notation $x \thickapprox y$ , if there is a $G$-equivariant
diffeomorphism $\Phi$ from an open $G$-invariant neighborhood $U$ of
$x$ onto an open $G$-invariant neighborhood $V$ of $y$. Clearly this
defines an equivalence relation $\thickapprox$ in $M$. The
equivalence classes will be called {\it orbit types} in $M$, and are
denoted by $M_x^\thickapprox$. If each stabilizer has only
finitely many components, then $x \thickapprox y$ if and only if
$G_x$ is conjugate to $G_y$ within $G$ and the actions of $G_x$, and
$G_y$, on $T_xM/T_xG(x)$, and $T_yM/T_yG(y)$, respectively, are
equivalent via a linear intertwining isomorphism (see chapter 2 of
\cite{DK}). The orbit $G(x)$ of $x\in M$ is {\it principal} if its type
$M_x^\thickapprox$ is open in $M$. Any non-principal orbit is called
a {\it singular} orbit. A nonprincipal orbit with the same dimension as a
principal orbit is an {\it exceptional} orbit.

Throughout the following $\R^3_1$ denotes the 3-dimensional real
vector space $\R^3$ with a scalar product of signature $(1,2)$ given
by  $\langle x,y\rangle=-x_1y_1 + x_2y_2+x_3y_3$.
 Let $Iso(\R^3_1)$ denote the group of isometries of $\R^3_1$, that is the Poincare
group $O(1,2)\ltimes \R^3$. We may write the natural action of an isometry
$(A,a) \in Iso(\R^3_1)$ as $(A,a)(x)=A(x)+a$,
where $A\in O(1,2)$ is called its {\it linear part} and $a\in \R^3$ is called its {\it translational part}.
Denote by $L:G\longrightarrow O(1,2)$
the projection on the linear part of $O(1,2)\ltimes \R^3$. If $L(G)$
is trivial then $G$ is called a {\it pure translation group}. We will
restrict our study to the identity component of $O(1,2)$, consisting
of orientation and time-orientation preserving isometries, which we denote
it by $SO_o(1,2)$, a subgroup of $O(1,2)$ of index 4.

In the next section we use an Iwasawa decomposition of $SO_o(1,2)$ to classify the Lie groups which act isometrically and by \co on $\R^3_1$. For this, we introduce a fixed Iwasawa decomposition of $SO_o(1,2)$. Let $i, j\in \{1,2,3\}$ and $E_{ij}$
be the $3\times 3$ matrix whose $(i,j)$-entry is $1$ and whose other entries are all $0$. Let $\fk= \{t(E_{23}-E_{32})|\ t\in \R\}$, $\mathfrak{a}=\{s(E_{12}+E_{21})|\ s\in \R\}$ and $\fn=\{u(E_{12}+E_{23}+E_{31}-E_{32})|\ u\in \R\}$. Then $\mathfrak{so{\rm (1,2)}}=\fk \oplus
   \mathfrak{a} \oplus \fn $  (direct sum  of vector spaces) is the Iwasawa decomposition of the Lie algebra $\mathfrak{so}{\rm (1,2)}$
   (see \cite[p. 372 ]{Kn}). Furthermore $SO_o(1,2)=KAN$ is the Iwasawa
    decomposition of $SO_o(1,2)$, in which $K$ , $A$ and $N$ are the
    connected Lie subgroups of $SO_o(1,2)$ associated to $\fk$ ,
    $\mathfrak{a}$ and $\fn$ respectively.
    Clearly, $\fk$ is isomorphic to
    $\mathfrak{so(2)}$ and using the exponential map, one gets that
    $K$ and $A$ are the standard embeddings of $SO(2)$ and $SO_o(1,1)$ in $SO_o(1,2)$, respectively. Each $C\in \fn$ is nilpotent, in fact $C^3=0$. It is
    easy to check that $[\mathfrak{a} , \fn ]\subseteq \fn$, so
    $\mathfrak{a} \oplus \fn$ is a Lie subalgebra of
    $\mathfrak{so{\rm (1,2)}}$ and $\fn$ is an ideal of $\mathfrak{a}
    \oplus \fn$. This implies that $N$ is normal in the group
    corresponding to $\mathfrak{a} \oplus \fn$, so $AN$ is a subgroup
    (in fact, nonabelian and solvable) of $SO_o(1,2)$. By a well known result any two dimensional Lie subgroup of $SO_o(1,2)$ is conjugate to $AN$.

 We end this section by introducing  he following notations.
Let $\{i, j\}\in\{1,2,3\}$. We fix the notations $B_i=E_{i(i+1)}+(-1)^{i+1}E_{(i+1)i}$, where $i\in\{1,2\}$, and $B_3=E_{13}+E_{31}+B_2$ throughout the paper. Let $\{e_1,e_2,e_3\}$ denote the standard basis of $\R^3$. Then $e_{ij}$ denotes the vector $e_i+e_j$, and $e_{123}:=e_1+e_2+e_3$.

\section{Lie groups acting by \co on $\R^3_1$}

In this section we determine all connected, closed Lie subgroups of $Iso_o(\R^3_1)$ acting isometrically and by \co on $\R^3_1$ up to conjugacy. If the Lie group $G$ is determined up to conjugacy, an immediate consequence is to specify the orbits
up to isometry. If $\dim L(G)=0$, since $L:G\rightarrow SO_o(1,2)$ is a Lie group
homomorphism and so continuous, $L(G)$ is connected and so it is trivial. In this case, $G$ is a two dimensional pure translation Lie subgroup of
$SO_\circ(1,2)\ltimes \R^3$ which its natural action on $\R^3_1$ is obviously, proper. If $\dim L(G)=1$, by a well-known result about one parameter Lie subgroups of $SO_o(1,2)$, the Lie group $L(G)$ is conjugate to one of the groups
$SO(2)\ (=K)$ , $A$ or $N$, where the representations of these groups in $SO_o(1,2)$ were introduced in the preceding section. We study the case, where $\dim L(G)=1$, in the following next three lemmas. In all of the lemmas, it is assumed that $G$ is a connected and closed Lie subgroup of $Iso(\R^3_1)$, which acts isometrically and by \co on $\R^3_1$.

\begin{lemma}\label{SO(2)} If  $L(G)$ is conjugate to $SO(2)$, then $G$ is conjugate to one of the following Lie groups within $SO_o(1,2)\ltimes \R^3$.

(i) The standard embedding of $SO(2)\times \R$ in $SO_o(1,2)\ltimes \R^3$.

(ii) The standard embedding of $Iso_o(\R^2)$ in $SO_o(1,2)\ltimes \R^3$.  \end{lemma}

{\bf Proof :} By the assumption $L(G)$ is conjugate to $SO(2)$. Then ,up to conjugacy, $l(\fg)=\{tB_2|t\in \R\}$. The action is of \co, so $\ker l$ should be one or two dimensional ideal of $\fg$. First assume that $\dim\ker l=1$. Then $\fg$, as a vector space, should be $\{(tB_2,(Du+D't)e_{123})|u,t\in \R\}$, where $D$ and $D'$ are two diagonal fixed $3\times 3$ matrices. Choose a vector $b$ in $\R^3$ such that $(B_2,b)\in \fg$. Since $\ker l$ is an ideal in $\fg$, we have $B_2(\ker l)\subseteq \ker l$. Hence $\ker l$ is a subspace of the eigenspace of $B_2$. The only one dimensional eigendirection of $B_2$ is $\R e_3$.  Hence $D_{22}=D_{33}=0$. Thus $(C,c)^{-1}\fg (C,c)=\{(tB_2,(D_{11}u+D'_{11}t)e_1)|u,t\in\R\}$, where $(C,c)=(I,D'_{33}e_2-D'_{22}e_3)\in SO_\circ (1,2)\ltimes \R^3$. Thus $G$ is conjugate to $SO(2)\times \R$.

Now suppose that $\dim\ker l=2$. This implies that $\dim G=3$. Then $\fg$, as a vector space, should be $\{tB_2,(Du+D'v+D''t)e_{123})|u,v,t\in \R\}$, where $D$, $D'$ and $D''$ are three diagonal fixed $3\times 3$ matrices. By rechoosing a basis for $\fg$ we may assume that $D_{33}=D'_{22}=0$, and up to conjugacy, $D''_{22}=D''_{33}=0$. We claim that both $D_{22}$ and $D'_{33}$ are nonzero. If $D'_{33}=0$. Choose two vectors $X_1=(E_{23}-E_{32},D''_{11}e_1)$ and $X_2=(0,(D+D')e_{123})$ from $\fg$. Then $[X_{1},X_2]=(0,-D_{22}e_3)$, which closeness under the bracket in $\fg$ implies that $D_{22}=0$. This implies that $\dim G=2$ which is in contradict to  $\dim G=3$. Hence $D'_{33}=0$. A similar discussion shows that $D_{22}\neq 0$. Hence without less of generality we may assume that $D_{22}=D'_{33}=1$. Now choose three vectors $X_1$ and $X_2$, as above, and $X_3=(0,(-A+B)e_{123})$ from $\fg$. By the fact that $[X_1,X_2]$ and $[X_1,X_3]$ belong to $\fg$, one gets that $D_{11}+D'_{22}=0$ and $D_{11}-D'_{22}=0$. So $D_{11}=D'_{22}=0$. If $D''_{11}\neq 0$, then $G(0)=\R^3$ which is in contradict to the assumption that the action is of cohomogeneity one. Hence $D''_{11}=0$. Thus $G$ is conjugate to the standard embedding of $Iso_o(\R^2)$ in $SO_o(1,2)\ltimes \R^3$.
{\it End of the proof of Lemma \ref{SO(2)}}.

 \begin{lemma}\label{A} If $L(G)$ is conjugate to $A$, then $\fg$ is conjugate to one of the following Lie algebras within $\mathfrak{so}(1,2)\oplus_\pi \R^3$, where $\pi:\mathfrak{so}(1,2)\rightarrow \mathfrak{Der}(\R^3_1)$ is the natural representation.

 $(i) \mathfrak{a}\oplus\R,$

 $(ii) \{(tB_1,ue_{1}+ve_2)|u,v\in \R\},$

$(iii) \{(tB_1,ue_{12}+ve_3|u,v\in\R\}$ or

$(iv)\{(tB_1, u(e_1-e_2)+ve_3|u,v\in\R\},$

$(v)\{(tB_1,s(e_1\pm e_2)+t\beta e_3)\ |\ t,s\in\R\},$\\
where $\beta$ is a fixed real number.
 \end{lemma}
   {\bf Proof :} By the assumption $l(\fg)=\{tB_1|t\in\R\}$, up to conjugacy. The action is of cohomopgeneity one, so $\ker l$ is a one or two dimensional ideal of $\fg$. First assume that $\dim\ker l=1$. Then by choosing a suitable coordinate we may assume that $\fg$, as a vector space, is $\{tB_1,(tD+sD')e_{123}\}$, where $D$ and $D'$ are two diagonal matrices. Now we determine the relation between the entries of $D$ and $D'$, to make $\fg$ as a Lie algebra.
    Take the following two vectors in $\fg$, $X_1=(B_1,De_{123})$ and $X_2=(0,D'e_{123})$. Then $[X_1,X_2]=(0,D'_{22}e_1+D'_{11}e_2)$, and so the closeness under the bracket implies the existence of a $s_0\in \R$ such that $D'_{11}=s_0D'_{22}$, $D'_{22}s_0=D'_{11}$ and $D'_{33}s_0=0$. Thus either $D'_{11}=D'_{22}=0$, $D'_{33}\neq 0$, or $D'_{11}=\pm D'_{22}\neq0$, $D'_{33}=0$. The first case shows that $(C,c)^{-1}G(C,c)=A\times \R$, where $(C,c)=(I,-(D_{22}e_1+D_{11}e_2))$, which is the case (i) of the lemma. The second case implies that $\fg$ is conjugate to one of the four Lie algebras $\{(tB_1,s(e_1\pm e_2)+tD_{33}e_3)\ |\ t,s\in\R\}$, depending on the real number $D_{33}$ to be zero or not, which is the case (v) of the lemma.

  Now suppose that $\dim\ker l=2$. Choose $b\in\R^3$ such that $(B_1,b)\in \fg$. Then $B_1(\ker l)\subseteq \ker l$. Hence $\ker l$ is a subspace of the eigenspace of $B_1$. This implies that $\ker l$ is one of the spaces $\{(0,ue_{1}+ve_2)|u,v\in \R\} ,\ \{(0,ue_{12}+ve_3|u,v\in\R\}\  {\rm or}\  \{(0, u(e_1-e_2)+ve_3|u,v\in\R\}.$ Thus $\fg$ is conjugate to one of the spaces $\{(tB_1,ue_{1}+ve_2)|u,v\in \R\}$,
$\{(tB_1,ue_{12}+ve_3|u,v\in\R\}$ or
$\{(tB_1, u(e_1-e_2)+ve_3|u,v\in\R\}$, which are the cases (ii) to (iv) of the lemma. {\it End of the proof of Lemma \ref{A}}.

\begin{lemma}\label{N} If $L(G)$ is conjugate to $N$, then $\fg$ is conjugate to one of the following Lie algebras within $\mathfrak{so}(1,2)\oplus_\pi \R^3$.

$(i) \{(tB_3,r(e_{12})+\beta te_3)| r,t\in\R\}$,

$(ii) \{(tB_3,r(e_{12})+se_3)|r,s,t\in \R\}$,\\
where $\beta$ is a fixed real number.

\end{lemma}
{\bf Proof :} By the assumption $l(\fg)=\{tB_3| t\in\R\}$, up to conjugacy. By a similar discussion of that of the proof of previous lemma, one gets that $\ker l$ is a one or two dimensional ideal of $\fg$. If $\dim\ker l=1$ then $\fg$, as a vector space, should be $\{(tB_3,(Du+D't)e_{123})|u,t\in \R\}$, where $D$ and $D'$ are two diagonal fixed matrices. By the fact that the only one dimensional eigendeirection of $B_3$ is $\R(e_{12})$ and using $B_3(\ker l)\subseteq \ker l$, one gets that $\fg=\{(tB_3,r(e_{12})+D'_{33}te_3)| r,t\in\R\}$, up to conjugacy. If $\dim\ker l=2$ then the relation $B_3(\ker l)\subseteq \ker l$ implies that $\ker l=\R (e_{12})+\R e_3$. Hence the \co assumption implies that $\fg=\{(tB_3,r(e_{12})+se_3)|r,s,t\in \R\}$. {\it End of the proof of Lemma \ref{N}}.

Any two dimensional Lie subgroup of $SO_o(1,2)$ is conjugate to $AN$. Hence for the case $\dim L(G)=2$ we have the following Lemma.
\begin{lemma}\label{AN} If $L(G)$ is conjugate to $AN$, then $\fg$ is conjugate to one of the following Lie algebras within $\mathfrak{so}(1,2)\oplus_\pi \R^3$.

$(i) (\mathfrak{a}\oplus \mathfrak{n})\times \{0\}$,

$(ii) \{(sB_1+tB_3,u(\alpha e_{12}+\beta e_3))|s,t,u\in\R\}$,

$(iii) \{(sB_1+tB_3,ue_{12}+ve_3))|s,t,u,v\in\R\}$,\\
where $\alpha$ and $\beta$ are two fixed real numbers.
\end{lemma}
{\bf Proof :} If $L(G)=AN$ then we have $$l(\fg)=\{sB_1+tB_3|s,t\in\R\}.$$
        By the \co assumption one gets that $0\leqslant\dim\ker l\leqslant 2$. If $\dim\ker l=0$, then $G$ is conjugate to $AN$ by a translation. If $\dim\ker l=1$, then $\ker l$, as a vector space, is $\{Due_{123}|u\in \R\}$, where $D$ is a fixed diagonal matrix. By the fact that $B_i(\ker l)\subseteq \ker l$, $i=1,3$, one gets that $D_{11}=D_{22}$. Hence, up to conjugacy, $$\fg=\{(sB_1+tB_3,u(D_{11}e_{12}+D_{33}e_3))|s,t,u\in\R\},$$ where $D_{11}$ and $D_{33}$ are fixed real numbers.

        If $\dim\ker l=2$, then $\ker l=\{(Du+D'v)e_{123}u,v\in \R\}$, where $D$ and $D'$ are two fixed diagonal matrices. By the same argument as above, one gets that $D_{11}=D'_{11}$ and $D_{22}=D'_{22}$. Hence by choosing a suitable basis for $\fg$ one gets that
        $$\fg=\{(sB_1+tB_3,ue_{12}+ve_3))|s,t,u,v\in\R\},$$
        up to conjugacy. {\it End of the proof of Lemma \ref{AN}}.

Since $\dim L(G)\leqslant 3$, the following lemma ends the classification of Lie groups acting isometrically and by \co on $\R^3_1$.

 \begin{lemma}\label{SO(1,2)} If $L(G)=SO_o(1,2)$, then $G=SO_o(1,2)\times \{0\}$.\end{lemma}
 {\bf Proof :}  If $L(G)=SO_o(1,2)$, then $\{B_1, B_2,B_3\}$ is a basis for $l(\fg)$. Since $B_i(\ker l)\subseteq \ker l$, where $i=1,2,3$, so $\ker l$ is either $\{0\}$ or $\R^3$. The action is of cohomogeneity one, so $\ker l=\{0\}$. Thus $G=SO_o(1,2)\times\{0\}$. {\it End of the proof of Lemma \ref{SO(1,2)}}.
\\

The main theorem of this section is the following.

\begin{theorem}\label{lg}
    Let $G$ be a closed and connected Lie subgroup of $Iso(\R^3_1)$ which acts isometrically and by \co on $\R^3_1$. Then the action is proper if and only if $G$ is conjugate to one of the following Lie groups.

    (a) A pure translation group,

    (b) The standard imbedding of $SO(2)\times \R$ in $SO_o(1,2)\ltimes \R^3$,

    (c) The standard imbedding of $Iso_o(\R^2)$ in $SO_o(1,2)\ltimes \R^3$,

    (d) $\{(A_t,u(e_1\pm e_2)+\beta te_3)\ |\ t,u\in\R\},$ where $\beta$ is a fixed nonzero real number.
\end{theorem}

   {\bf Proof} : By Lemmas \ref{SO(2)} to \ref{SO(1,2)} we know the Lie algebras of all Lie groups which act isometrically and by \co on $\R^3_1$.
   Hence, to prove of the theorem, we need only to investigate those acting properly. If $\dim L(G)=0$ then $G$ is a pure translation group and the action is obviously proper. If $L(G)=SO(2)$ then $G$ is conjugate to either $SO(2)\times \R$ or $Iso_o(\R^2)$ and so the action is reduced to the action of a closed Lie subgroup of $Iso(\mathbb{E}^3)$, which is proper clearly. We claim that if $L(G)$ is noncompact then case (d) of the theorem occur. If $L(G)$ is not compact then $\fg$ is conjugate to one of the Lie algebras stated in Lemmas \ref{A} to \ref{SO(1,2)}. All the Lie algebras listed in Lemma \ref{A}-(i) to (iv), Lemma \ref{N}-(ii), Lemma \ref{AN} to \ref{SO(1,2)} cause a nonproper action, since in each case the stabilizer of the origin is not compact. If $\fg$ is conjugate to that of in Lemma \ref{N}-(i), then $\exp (\fh)$ is a closed and noncompact subgroup of the stabilizer of each point $(x,x+\beta,0)^{T}\in\R^3_1$, where $\fh=\{(tB_3,\beta te_3)| t\in\R\}$. This shows that the action is nonproper in case of Lemma \ref{N}-(i). Thus to complete the proof of our claim, we only need to verify that the action caused by the Lie algebra stated in Lemma \ref{A}-(v) is proper, if $\beta$ is nonzero. In that case, by a simple computation one gets that $$G=\exp (\fg)=\{(A_t,u(e_1\pm e_2)+\beta te_3)\ |\ t,u\in\R\}.$$
   Let $\{t_n\}$ and $\{u_n\}$ be two real sequences. Let $\{X_n=(x_n,y_n,z_n)^T\}$ be a sequence in $\R^3_1$. Let $\{g_n=(\cosh t_n(E_{11}+E_{22})+\sinh t_n(E_{12}+E_{21})+E_{33},u_n(e_1\pm e_2)+\beta t_ne_3)\}$ be a sequence in $G$. Let $g_n.X_n\rightarrow Y$ and $X_n\rightarrow X$, when $n\rightarrow +\infty$. If $Y=(y_1,y_2,y_3)^T$ and $X=(x_1,x_2,x_3)^T$, then $t_n\rightarrow \frac{y_3-x_3}{\beta}$ and
$u_n\rightarrow y_1-x_1\cosh \frac{y_3-x_3}{\beta}-x_2\sinh \frac{y_3-x_3}{\beta}$. Hence $\{g_n\}$ is a convergent sequence in $G$, and so the action is proper. {\it End of the proof of Theorem \ref{lg}}.

As a consequence of Lemmas \ref{A} to \ref{AN} and Theorem \ref{lg} one gets the following corollary. The Lie groups in the list are obtained by the exponential function.
\begin{cor}\label{nonproper}
 Let $G$ be a closed and connected Lie subgroup of $Iso(\R^3_1)$ which acts isometrically and by \co on $\R^3_1$. Then the action is nonproper if and only if $G$ is conjugate to one of the following Lie groups within $SO_o(1,2)\ltimes \R^3$.

 $(i)\ \{(A_t,se_3)|t,s\in \R \}$,

$(ii)\ Iso_o(\R^2_1)$,

$(iii)\ \{(A_t,ue_{12}+ve_3)|t,u,v\in \R \}$,

$(iv)\ \{A_t, u(e_1-e_2)+ve_3|t,u,v\in\R\}$,

$(v)\ \{(A_t,ue_{12})|t,u\in \R \}$,

$(vi)\ \{(A_t,u(e_1-e_2)|t,u\in \R \}$,

$(vii)\ \{(N_t,ue_{12}+\beta te_3)|t,u\in \R \}$,

$(viii)\ \{(N_t,ue_{12}+ ve_3)|t,u,v\in \R \}$,

$(ix)\ AN$,

$(x)\ AN\ltimes_\pi\{u(\alpha e_{12}+\beta e_3)|\ u\in \R\}$,

$(xi)\ AN\ltimes_\pi\{u e_{12}+ve_3)|\ u,v\in \R\}$,

$(xii)\ SO_o(1,2)$,\\
where $\alpha$ and $\beta$ are fixed real numbers and $\pi: AN\rightarrow \R^3$ is the natural representation.
\end{cor}
 \section{Causal properties of the orbits}
Assume that the connected and closed Lie subgroup $G$ of
$Iso(\R^3_1)$ acts isometrically and by \co on $\R^3_1$, we determine
causal properties of the orbits.

The orbit $G(p)$ is said to be Lorentzian, degenerate or space-like
if the induced metric on $G(p)$ is Lorentzian, degenerate or
Riemannian, respectively. It is called time-like or light-like if
each nonzero tangent vector in $T_pG(p)$ is time-like or null,
respectively. The category into which a given orbit falls is called its {\it causal property}.
\subsection{The action is proper}
Let a Lie group $G$ act by \co and properly on a smooth manifold
$M$.  A result by Mostert (see \cite {Mos}), for the compact Lie
groups, and Berard Bergery (see \cite {Ber}), for the general case,
says that the orbit space $M/G$ is homeomorphic to one of the spaces
$$\R \quad, \quad S^1 \quad, \quad [0,+\infty) \quad , \quad
[0,1].$$ In the following theorem we show that the cases $[0,1]$ and $S^1$ can
not occur, when $M=\R^3_1$.

\begin{theorem}\label{orbit}
 Let $\R^3_1 $ be of \co under the isometric action of a connected
      and closed Lie subgroup $G\subset Iso(\R^3_1)$. If the action is proper, then
      one of the following cases occurs.

      (1) $G$ is a pure translation group. In this case each orbit is a plane which is obtained by a translation of $G(0)$, and so the orbit space is diffeomorphic to $\R$.

      (2) $G$ is conjugate to $SO(2)\times \R$. In this case, there is a time-like singular orbit which is a one dimensional affine subspace and each other orbit is a cylinder about the singular orbit and so the orbit space is diffeomorphic to $[0,+\infty)$. In particular, each principal orbit is a Lorentzian cylinder.

      (3) $G$ is conjugate to $Iso_o(\R^2)$. In this case, each orbit is a space-like plane which is obtained by a translation of $G(0)$, and so the orbit space is diffeomorphic to $\R$.

      (4) $G$ is conjugate to   $\{((E_{11}+E_{22})\cosh t+(E_{12}+E_{21})\sinh t+E_{33},u(e_1\pm e_2)+\beta te_3)\ |\ t,u\in\R\}$,
where $\beta$ is a fixed nonzero real number. In this case an orbit is a degenerate plane, and each other orbit is a Lorentzian generalized cylinder. The orbit space is diffeomorphic to $\R$.
\end{theorem}
     {\bf Proof} : The theorem is a direct consequence of Theorem \ref{lg}, and only the case (4) needs some explanation. Suppose that $\fg=\{(t(E_{12}+E_{21}),s(e_{12})+t\beta e_3)\ |\ t,s\in\R\}$. Let $X_1=(E_{12}+E_{21},\beta e_3)$ and $X_2=(0,e_1+ e_2)$. Then $\{X_1,X_2\}$ is a basis for $\fg$. Fix an arbitrary point $X=(x,y,z)^T \in \R^3_1$, then
          $$\frac{d}{dt}\exp(tX_1)(X)|_{t=0}=(x\sinh t
          +(y+\beta)\cosh t)e_1 + (x\cosh t +(y+\beta)\sinh t)e_2 + e_3, $$
and

          $\frac{d}{ds}\exp(sX_2)(X)|_{s=0}=e_{12}.$

          If $x=y$ then the vector $N=e_1+e_2$ is normal to the above two vectors, and
          so to $G(X)$ which implies that $G(X)$ is a degenerate principal orbit. It is easily seen that this orbit is a plane.

          If $x\neq y$ then the unit space-like vector $N=\frac{e^{t}}{y+\beta -x}(1,1,(y+\beta -x)e^{-t})^T$ is
          normal to $G(X)$ at $X$, so $G(X)$ is a Lorentzian orbit,
          i.e. $T_XG(X)$ is isometric to $\R^2_1$,
          and the shape operator associated to $N$ is represented
          with respect to the pseudo-orthogonal basis
          $v_1=\sqrt{2}/2(1,1)^T$ , $v_2=\sqrt{2}/2(1,-1)^T$ as follows
          $$S=\left[\begin{array}{cc} 0 & 0\\ 1 & 0
          \end{array}\right].$$
          Hence the shape operator is not diagonalizable. Thus $G(X)$ is locally isometric
          to a generalized cylinder by \cite {M}. {\it End of the proof of Theorem \ref{orbit}}.

\begin{cor}\label{orbit2}
 Let $\R^3_1 $ be of \co under the proper action of a connected
      and closed Lie subgroup $G\subset Iso(\R^3_1)$.  then one of the followings holds:

       1) If there is a singular orbit, then it is a time-like one dimensional
     affine subspace and each principal orbit is isometric to $\R^1_1\times S^1(r)$ for some
      $r>0$.

      2) If there is a space-like orbit, then each orbit is a space-like hyperplane.

      3) If there is a Lorentzian orbit isometric to $\R^2_1$, then
      each orbit is isometric to $\R^2_1$.

      4) If there are more than one degenerate orbit, then each
      orbit is a degenerate hyperplane.

      5) If there is exactly one degenerate orbit, then it is a degenerate
      hyperplane and each other orbit is locally isometric to a generalized cylinder and there is no singular
       orbit.
\end{cor}

\subsection{The action is nonproper}
As a consequence of Corollary \ref{nonproper} one gets that if the action is not proper then $L(G)$ is conjugate to one of the Lie groups $A$, $N$, $AN$ or $SO_o(1,2)$. In the following propositions we consider the action of each of the Lie groups, mentioned in Corollary \ref{nonproper}, and then investigate the causal properties of the orbits. Throughout this subsection $p=(x_0,y_0,z_0)^T$ denotes an arbitrary fixed point in $\R^3_1$.
\begin{pro}\label{pro-nonproper-1}
Let $\R^3_1 $ be of \co under the isometric and nonproper action of a connected
and closed Lie subgroup $G\subset Iso(\R^3_1)$. If $L(G)=A$, then one of the following cases occurs.

(i) There is a one dimensional space-like affine subspace as a singular orbit. There are four degenerate half-plans which are principal orbits. Each other orbit is a branch of a Lorentzian hyperbolic cylinder which is principal.

(ii) Each orbit is a translation of a Lorentzian plane.

(iii) There is a degenerate plane as a unique exceptional orbit and there are two open submanifolds as the orbits.

(iv) There are infinitely many one dimensional light-like singular orbits. Each other orbit is a Lorentzian principal orbit, which is not closed submanifold.
\end{pro}

{\bf Proof} :  By Corollary \ref{nonproper} $G$ is conjugate to one of the stated Lie groups in $(i)$ to $(vi)$.
(i) Let $G=\{(A_t,se_3)|t,s\in \R \}$. If $x_0=y_0=0$, then $G(p)=\{se_3|s\in \R\}$ and so $G(p)$  is a one dimensional space-like singular orbit. Let $x\neq 0$ or $y\neq 0$. Then $G(p)$ is a principal orbit, since $G_p$ is the trivial subgroup. It is easily seen that if $(u,v,w)^T$  belongs to $G(p)$, then $u^2-v^2=x_0^2-y_0^2$. Hence there are four degenerate principal orbits for the cases $x_0=\pm y_0$ and depending on the sign of $x_0$. If $x_0\neq \pm y_0$, then $G(p)$ is a branch of a Lorentzian hyperbolic cylinder $u^2-v^2=x_0^2-y_0^2$.

(ii) Let $G=Iso_o(\R^2_1)$. It is easily seen that each orbit is a translation of the Lorentzian plane $z=0$.

(iii) Let $G=\{(A_t,ue_{12}+ve_3)|t,u,v\in \R \}$. If $x_0=y_0$ then $G(p)=\{ue_{12}+ve_3|\ u,v\in\R\}$ and $G_p=\{(A_t,x(e^t-1)e_{12})|\ t\in\R\}$, and so $G(p)$ is a degenerate plane as an exceptional orbit.  If $x_0\neq y_0$, then $G_p$ is the trivial subgroup and so $G(p)$ is an open submanifold. So there are two open orbits corresponding to the cases $x_0>y_0$ and $x_0<y_0$. A similar discussion about the action of the Lie group $G=\{A_t, u(e_1-e_2)+ve_3|t,u,v\in\R\}$, stated in Corollary \ref{nonproper}-(iv), yields the same results stated in Proposition \ref{pro-nonproper-1}-(iii).
%

(iv) Let $G=\{(A_t,ue_{12})|t,u,v\in \R \}$. If $x_0=y_0$ then $G(p)=\{ue_{12}+z_0e_3|\ u\in\R\}$. Hence $(x,y,z)^T$ belongs to $G(p)$ if and only if $x=y$ and $z=z_0$. Hence there are infinitely many one dimensional light-like singular orbits. If $x_0\neq y_0$, then $G_p$ is the trivial subgroup and so $G(p)=\{ue_1+ve_2+z_0e_3|\ (u-v).(x_0-y_0)>0\}$. This shows that $G(p)$ is a Lorentzian principal orbit, which is not closed submanifold. A similar discussion about the action of the Lie group $G=\{(A_t,u(e_1-e_2)|t,u\in \R \}$, stated in Corollary \ref{nonproper}-(vi), yields the same results stated in Proposition \ref{pro-nonproper-1}-(iv). {\it End of the proof of Proposition \ref{pro-nonproper-1}.}

\begin{pro}\label{pro-nonproper-2}
Let $\R^3_1 $ be of \co under the isometric and nonproper action of a connected
and closed Lie subgroup $G\subset Iso(\R^3_1)$. If $L(G)=N$, then one of the following cases occurs.

(i) The acting group is conjugate to $\{(N_t,ue_{12}+\beta te_3)|t,u\in \R \}$. There are two cases. If $\beta\neq 0$ then each orbit is principal which is obtained by a translation of a fix degenerate plane. If $\beta=0$ then there are infinitely many light-like one dimensional singular orbits of the same type, where the union of them is a degenerate plane. Each principal orbit is obtained by a translation of a fix degenerate plane. In both cases all the principal orbits are of the same type.

(ii) The acting group is conjugate to $\{(N_t,ue_{12}+ ve_3)|t,u,v\in \R \}$. Each orbit is a degenerate plane as a principal orbit and the set of the orbits is a foliation of $\R^3_1$. All of the orbits are of the same type.
\end{pro}
{\bf Proof:} By Corollary \ref{nonproper} $G$ is conjugate to one of the stated Lie groups in $(vii)$ and $(viii)$. So we have the following two cases.

(i) Let $G=\{(N_t,ue_{12}+\beta te_3)|t,u\in \R \}$. If $x_0\neq y_0$ then $G_p$ is trivial and $G(p)=p+\{ue_{12}+ve_3|\ u,v\in\R\}$, which is a degenerate plane. Since the union of these orbits is an open subset of $\R^3_1$, so each of them is principal. Let $x_0=y_0$. If $\beta=0$, then $G_p=\{(N_t,(-tz_0)e_{12})|\ t\in \R\}$. This shows that $G(p)$, which is equal to $p+\{ue_{12}|\ u\in\R\}$, is a one dimensional light-like subspace as a singular orbit. Let $p'=(x,y,z)^T$. It is easily seen that $p'$ belongs to $G(p)$ if and only if $x=y$ and $z=z_0$. Furthermore, if $p'\notin G(p)$ and $x=y$, then $G_{p'}=g^{-1}G_{p}g$, where $g=(I,(z+z_0)e_3)$. Thus  there are infinitely many singular orbits of the same type. Obviously, the union of the singular orbits is $\{ue_{12}+ve_3|\ u,v\in\R\}$, which is a degenerate plane. If $\beta\neq0$ (and $x_0=y_0$), then $G_p$ is the trivial subgroup and so $G(p)=2$ is a degenerate plane as a principal orbit.

(ii) Let $G=\{(N_t,ue_{12}+ ve_3)|t,u,v\in \R \}$. Then $G_p=\{N_t,((y_0-x_0)\frac{t^2}{2}-tz_0)e_{12}+(y_0-x_0)te_3\}$, and $G(p)=p+\{ue_{12}+ ve_3|\ u,v\in\R\}$. Hence each orbit is a degenerate plane as a principal orbit. All of the orbits are of the same type, since each of the stabilizers is conjugate to $\{(N_t,0)|\ t\in\R\}$. In fact $g^{-1}G_pg=\{(N_t,0)|\ t\in\R\}$, where $g=(I,(y_0-x_0)e_2+z_0e_{3})$. {\it End of the proof of Proposition \ref{pro-nonproper-2}.}
\begin{pro}\label{pro-nonproper-3}
Let $\R^3_1 $ be of \co under the isometric and nonproper action of a connected
and closed Lie subgroup $G\subset Iso(\R^3_1)$. If $L(G)=AN$, then one of the following cases occurs.

(i) The acting group is conjugate to $AN$. Then there is one principal orbit type and two singular orbit types. Each principal orbit is either a Lorentzian or a space-like surface. There is a zero dimensional singular orbit and infinitely many one dimensional light-like singular orbits.

(ii) The acting group is conjugate to $G=AN\ltimes_\pi\{u(\alpha e_{12}+\beta e_3)|\ u\in \R\}$. If $\alpha =0$, then there is a degenerate exceptional orbit and there are two orbits which are open submanifolds of $\R^3_1$ (so the orbit space consists of three points). If $\alpha\neq 0$, then there are infinitely many one dimensional light-like singular orbits of the same type and two orbits which are open submanifolds of $\R^3_1$.

(iii) The acting group is conjugate to $G=AN\ltimes_\pi\{u e_{12}+ ve_3|\ u,v \in \R\}$. Then there is a degenerate exceptional orbit and there are two orbits which are open submanifolds of $\R^3_1$ (so the orbit space consists of three points).
\end{pro}
{\bf Proof :} By Corollary \ref{nonproper} $G$ is conjugate to one of the stated Lie groups in $(ix)$ to $(xi)$. So we have the following two cases.

(i) Let $G=AN$. Then $G$ fixes the origin, and so the origin is a singular orbit. Let $p$ be not the origin. The set $\{B_1,B_3\}$ is a basis for the Lie algebra $\fg$. To determine causal properties of the orbits, let
$$\Phi_p(t)=\exp ((t\alpha B_1+B_3)).p,$$
where $\alpha$ is an arbitrary real number. Then
\begin{equation}\label{equ-AN}
\langle \frac{d\Phi_p}{dt}(0), \frac{d\Phi_p}{dt}(0)\rangle =\alpha^2 (x_0^2-y_0^2)+2\alpha z_0(x_0-y_0)+(x_0-y_0)^2.
\end{equation}
This implies that if $x_0\neq y_0$ and $p$, as a vector, is a nonzero space-like (resp. time-like) vector, then the polynomial (\ref{equ-AN}) has two roots (resp. has no root). Hence the orbit $G(p)$ is a Lorentzian (resp. space-like) orbit. Since $G_p$ is trivial and the union of these orbits is an open subset of $\R^3_1$, all of these orbits are of the same type and principal. If $x_0=y_0$, then $G(p)\{ue_{12}+z_0e_3|\ ux_0>0\}$, and so it is a one dimensional light-like singular orbit. All of these orbits are of the same type, since their stabilizers are conjugate to $N$.

(ii) Let $G=AN\ltimes_\pi\{u(\alpha e_{12}+\beta e_3)|\ u\in \R\}$.

Let $\alpha =0$. If $x_0=y_0$ then $G(p)=\{ue_{12}+ve_3|\ u,v\in \R\}$, which is a two dimensional degenerate orbit. If $x_0\neq y_0$ then $\dim G(p)=3$ and $(u,v,w)^T\in G(p)$ if and only if $(u-v)(x_0-y_0)>0$. Hence there are two open submanifolds as the three dimensional orbits. This implies that the two dimensional degenerate orbit is an exceptional orbit.

 Let $\beta =0$ (and so $\alpha\neq 0$). If $x_0=y_0$ then $G(p)=p+\{ue_{12}|\ u\in \R\}$, which is a one dimensional light-like singular orbit. Furthermore, $(u,v,w)^T\in G((x,y,z)^T)$ if and only if $w=z$. This implies that there are infinitely many singular orbits. Since their stabilizers are conjugate to $AN$, all of them are of the same type. If $x_0\neq y_0$, then $\dim G(p)=3$, and $(u,v,w)^T\in G(p)$ if and only if $(u-v)(x_0-y_0)>0$. Hence there are two open submanifolds as the three dimensional orbits.

 Now let both of $\alpha$ and $\beta$ are nonzero. If $x_0=y_0$ then $G(p)=p+\{\alpha ue_{12}+\beta ue_3|\ u\in \R\}$, which is a one dimensional light-like singular orbit. A similar discussion of that of the previous case, shows that there are infinitely many singular orbits of the same type and two open orbits in $\R^3_1$.

 (iii) Let $G=AN\ltimes_\pi\{u e_{12}+ ve_3|\ u,v \in \R\}$.  If $x_0=y_0$ then $G(p)=\{ue_{12}+ ve_3|\ u,v \in \R\}$. If $x_0\neq y_0$ then $\dim G(p)=3$.  In the later case, $(u,v,w)^T\in G((x,y,z)^T)$ if and only if $(u-v)(x-y)>0$. Hence there is a degenerate exceptional orbit and two orbits which are open submanifolds of $\R^3_1$. {\it End of the proof of Proposition \ref{pro-nonproper-3}}.
\\

By Corollary \ref{nonproper} the only case, that we have not investigated in three previous propositions, is the case that $G$ is conjugate to $SO_o(1,2)$. Let $G=SO_o(1,2)$ and let $G$ act on $\R^3_1$ naturally. Then the origin is a zero dimensional singular orbit, each component of the light-cone is an exceptional orbit, and each pseudo-sphere and each pseudo-hyperbolic space is a principal orbit. Hence there is one singular orbit type, one exceptional orbit type and two principal orbit types. By reviewing Propositions \ref{pro-nonproper-1} to \ref{pro-nonproper-3} we get that there are at most two exceptional orbits and we obtain the following corollary.
\begin{cor}
Let $\R^3_1 $ be of \co under the isometric and nonproper action of a connected
and closed Lie subgroup $G\subset Iso(\R^3_1)$. If there is a unique exceptional orbit then it is a degenerate plane and there are two orbits which are open submanifolds. In particular the orbit space, which consists of three points, is not Housdorff.
\end{cor}


 \bibliographystyle{amsplain}

\end{document}